\documentclass[10pt, conference]{IEEEtran}
\usepackage{epsfig,color,amsmath,cite}
\usepackage{amsthm} 
\usepackage{amsmath}    
\IEEEoverridecommandlockouts
\usepackage{bm}
\usepackage{epstopdf}
\usepackage{amssymb}
\usepackage{url}
\usepackage{enumitem} 
\usepackage{multirow}
\usepackage{hhline}
\usepackage{booktabs}
\usepackage{mathtools}
\usepackage[linesnumbered,boxed,commentsnumbered,ruled,vlined,longend]{algorithm2e}
\usepackage{comment}

\DeclareMathOperator*{\minimize}{minimize}

\DeclareMathOperator*{\subjectto}{subject\ to}

\makeatother
\DeclareMathAlphabet\mathbfcal{OMS}{cmsy}{b}{n}

\newtheorem{theorem}{Theorem}

\newtheorem{asmp}{Assumption}



\usepackage{stackengine}
\newcommand\barbelow[1]{\stackunder[1.2pt]{$#1$}{\rule{.8ex}{.075ex}}}

\newcommand{\mat}[1]{\boldsymbol{#1}}

\newcommand{\bmat}[1]{\begin{bmatrix} #1 \end{bmatrix}}

\providecommand{\mA}{\ensuremath{\mat{A}}}
\providecommand{\mB}{\ensuremath{\mat{B}}}
\providecommand{\mC}{\ensuremath{\mat{C}}}

\providecommand{\mL}{\ensuremath{\mat{L}}}

\providecommand{\mP}{\ensuremath{\mat{P}}}

\providecommand{\mY}{\ensuremath{\mat{Y}}}

\newcommand{\m}{\boldsymbol}
\allowdisplaybreaks[4]
\usepackage[colorlinks = true,
linkcolor = blue,
urlcolor  = blue,
citecolor = blue,
anchorcolor = blue]{hyperref}

\newcommand{\mbb}[1]{\mathbb{#1}}

\usepackage[framemethod=TikZ]{mdframed}
\mdfdefinestyle{MyFrame}{%
	linecolor=black,
	outerlinewidth=1.25pt,
	roundcorner=1.25pt,
	innerrightmargin=5pt,
	innerleftmargin=5pt,}
	
\usepackage[noabbrev]{cleveref}

\usepackage{mathtools}

\DeclarePairedDelimiter\abs{\lvert}{\rvert}%
\DeclarePairedDelimiter\norm{\lVert}{\rVert}%

\makeatletter
\let\oldabs\abs
\def\abs{\@ifstar{\oldabs}{\oldabs*}}
\let\oldnorm\norm
\def\norm{\@ifstar{\oldnorm}{\oldnorm*}}
\makeatother

\usepackage[english]{babel}
\usepackage[utf8]{inputenc}
\usepackage[super]{nth}

\usepackage{graphicx}
\usepackage{float}
\usepackage[caption = false]{subfig}

\usepackage{array}
\usepackage{threeparttable}

\usepackage[english]{babel}
\usepackage[utf8]{inputenc}
\usepackage[super]{nth}

\RequirePackage{filecontents}

\usepackage{lipsum}

\title{\vspace{0.7cm} \centering \LARGE \textbf{Sensor Placement Strategies for Some Classes of Nonlinear \\ Dynamic Systems via Lyapunov Theory}}

\author{Sebastian A. Nugroh$\text{o}^{\star}$ and Ahmad F. Tah$\text{a}^{\star}$ \vspace{-0.25cm}
	\thanks{
		$^{\star}$Department of Electrical and Computer Engineering, The University of Texas at San Antonio, 1 UTSA Circle, San Antonio, TX 78249.
		Emails: sebastian.nugroho@my.utsa.edu, ahmad.taha@utsa.edu. This work is partially supported by the National Science Foundation under Grants 1728629 and 1917164.}
}

\begin{document}

\newdimen\origiwspc%
\newdimen\origiwstr%
\origiwspc=\fontdimen2\font
\origiwstr=\fontdimen3\font

\fontdimen2\font=0.63ex

\maketitle

\begin{abstract}
In this paper, the problem of placing sensors for some classes of nonlinear dynamic systems (NDS) is investigated. In conjunction with mixed-integer programming, classical Lyapunov-based arguments are used to find the minimal sensor configuration such that the NDS internal states can be observed while still optimizing some estimation metrics. The paper's approach is based on two phases. The first phase assumes that the encompassed nonlinearities belong to one of the following function set classifications: \textit{bounded Jacobian}, \textit{Lipschitz continuous}, \textit{one-sided Lipschitz}, or \textit{quadratically inner-bounded}. 
To parameterize these classifications, two approaches based on stochastic point-based and interval-based optimization methods are explored. Given the parameterization, the second phase formulates the sensor placement problem for various NDS classes through mixed-integer convex programming. The theoretical optimality of the sensor placement alongside a state estimator design are then given. Numerical tests on traffic network models showcase that the proposed approach yields sensor placements that are consistent with conventional wisdom in traffic theory.
\end{abstract}

%
%
%

\section{Introduction and Motivation}

The selection/placement problem of sensors and control nodes (SC) in dynamic networks has received noticeable attention in the past decades from a variety of research communities spanning systems biology, physics, social sciences, and engineering.
The optimal selection/placement of SC in {linear} networks has been thoroughly studied in the literature through three major approaches.
The first approach is based on combinatorial heuristics and detailed routines that exploit network structure and properties~\cite{tzoumas2016minimal,zhang2017sensor,pequito2016}. The second approach entails utilizing semidefinite programming (SDP) formulations of control/estimation methods while including sparsity promoting penalties---thereby minimizing the total number of activated SC~\cite{Dhingra2014,Argha2017}. The third approach uses mixed-integer convex programming, convex approximations, and convex relaxations to obtain the minimal set of SC~\cite{chanekar2017optimal,Taha2018,Singh2018,chang2018co,Ebrahimi2018}. 

As for \textit{nonlinear} systems, the  \textit{sensor placement problem} (SPP) remains an open
 research problem. Solving the SPP for a nonlinear representation of a dynamic network---in contrast with linearizing around an operating point---can offer a sensor selection that \textit{works} for all operating regions of the network. Our work-in-progress abstract \cite{Nugroho2019sap} showcases this particular advantage: solving the SPP based on a nonlinear representation of stepper motor offers a linearization-free placement that is feasible for any operating points. Strategies to solve the SPP in NDS that are available in literature can be summarized as follows. First, Qi \textit{et al.} \cite{qi2015optimal} utilize the concept of \textit{empirical} Gramians to obtain a measure towards observability in NDS, which is then used to determine the location of {phasor measurement unit} in power networks. Next, Haber \textit{et al.}\cite{Haber2017} present a method for reconstructing the initial states $\m x(t_0)$ of NDS while optimally selecting sensors for a given observation window. 
It has been argued in \cite{Haber2017} that the proposed approach is computationally more tractable than the one in~\cite{qi2015optimal}.
Recently, a randomized algorithm for dealing with SPP is introduced by Bopardikar \textit{et al.}\cite{Bopardikar2019}, where the authors develop theoretical bounds for eigenvalue and condition number of observability Gramian. 
Still, the potential applicability of the aforementioned methods to address the SPP on large-scale, unstable NDS along with their feasibility to incorporate sensor selection with some estimation metrics all remain unclear.

To that end, in this paper we focus on addressing the SPP for NDS using a more general framework by firstly incorporating various observer designs---typically designed using Lyapunov theory---available in the literature. Secondly, we exploit two important \textit{attributes} that are prevalent in the majority of NDS: \textit{(i)} physical states in many dynamic networks are almost always bounded (e.g., traffic in transportation networks, tank levels and flows in water networks, frequencies in power grids), and \textit{(ii)} the nonlinearities of NDS can be shown---or rather assumed---to belong to some classes of nonlinear function sets.
In what follows we summarize the paper's approach and contributions.


Section~\ref{sec:problem_formulation} formalizes the SPP for NDS. In Section~\ref{sec:NDS_classification}, by exploiting the aforementioned attributes, we discuss two distinct methods to parameterize the nonlinearities in NDS 
that can be classified as \textit{(a)} \textit{bounded Jacobian}, \textit{(b)} \textit{Lipschitz continuous}, \textit{(c)} \textit{one-sided Lipschitz}, and \textit{(d)} \textit{quadratically inner-bounded}. 
The parameterization---performed via stochastic point-based optimization and interval-based global optimization---essentially finds constants that bound the nonlinearities. Section~\ref{sec:sensor_placement}, and given the parameterization, formulates the SPP through mixed-integer SDP (MISDP) for Lipschitz NDS, which gives a minimal sensor placement for Lipschitz nonlinear systems together with a Luenberger-like observer gain. Since we use a generalized framework, as shown later in Section~\ref{sec:sensor_placement}, this approach can be easily extended to formulate the SPP for other classes of NDS as well as other type of observer designs.
Section~\ref{sec:numerical} illustrates the effectiveness of the proposed approach to address the SPP on traffic networks. The summary of this paper along with potential future work are all provided in Section~\ref{sec:conclusion_future_work}. 

\noindent \textit{\textbf{Notation}}: Italicized, boldface upper and lower case characters represent matrices and column vectors: $a$ is a scalar, $\m a$ is a vector, and $\m A$ is a matrix. Matrix $\m I$ denotes the identity square matrix---matrix $\m I_{n}$ specifies the size of identity matrix of dimension $n\times n$. The notations $\m a_n$ and $\m a_{n\times m}$ denote $n\times 1$ vector and $n\times m$ matrix with all elements are equal to $a$, while $\mathbb{R}$ and $\mathbb{R}_+$ denote the set of real and non-negative real numbers. The notations $\mathbb{R}^n$ and $\mathbb{R}^{p\times q}$ denote the sets of row vectors with $n$ elements and matrices with size $p$-by-$q$ with elements in $\mathbb{R}$. For any $\m x \in \mathbb{R}^{n}$, $\Vert\m x\Vert_2$ denotes the Euclidean norm of of $\m x$, where $\m x^{\top}$ denotes the transpose of $\m x$. If $\m A$ is a matrix, its $i$-th and $j$-th element is denoted by $A_{(i,j)}$. The operators $\mathrm{blkdiag}(\cdot)$ constructs a block diagonal matrix, $\mathrm{diag}(\cdot)$ constructs a diagonal matrix from a vector, $\mathrm{vec}(\cdot)$ constructs a vector by stacking each column of a matrix, $\otimes$ denotes the Kronecker product, and $\oslash$ denotes the Hadamard product. 

\section{The Sensor Placement Problem (SPP) for NDS}\label{sec:problem_formulation}
We consider NDS modeled in the following form
\begin{subequations}\label{eq:gen_dynamic_systems}
\begin{align}
\dot{\m x}(t) &= \mA \m x (t) + \m f(\m x) + \mB\m u(t),\\
\m y(t) &= \mC \m x (t),
\end{align}
\end{subequations}
\noindent where vectors $\m x\in \mathbfcal{X}$, $\m u\in\mathbb{R}^{n_u}$, $\m y \in \mathbb{R}^{n_y}$ respectively denote the state, input, and output of the system, function $\m f :\mathbb{R}^{n_x}\rightarrow \mathbb{R}^{n_x}$ represents any existing nonlinearity, and $\m A$, $\mB$, $\mC$ are constant matrices of appropriate dimensions. In the sequel we occasionally drop the time dependence such that $\m x := \m x(t)$ to save space. 
We assume that there are $N$ nodes where each node $i$ has $n_{x_i}$ number of states, $n_{u_i}$ number of inputs, and $n_{y_i}$ number of outputs. For simplicity, we assume that the mapping from state to output vector $\m C \m x(t)$ in \eqref{eq:gen_dynamic_systems} has the following structure $\m C \m x = \mathrm{blkdiag}\left(\m C_1,\hdots,\m C_N\right)\left[\m x_1^{\top}\;\cdots\;\m x_N^{\top}\right]^{\top}$.
This structure allows us to conveniently select particular nodes which output are measured. To construct the SPP for NDS  \eqref{eq:gen_dynamic_systems}, we consider a binary variable, denoted by $\gamma_i$, that determines the activation or deactivation of sensor on each node $i$---we set $\gamma_i = 1$ if the sensor measuring node $i$ is activated and $\gamma_i = 0$ otherwise. These variables can be compactly organized into one vector namely $\m \gamma$ such that $\m \gamma = \left[\gamma_1\,\cdots\,\gamma_N\right]^{\top}$. In addition, we also define $\m \Gamma(\m \gamma) := \mathrm{blkdiag}\big(\gamma_1\m I_{n_{y_1}},\hdots,\gamma_N\m I_{n_{y_N}}\big)$ such that the NDS along with sensor placement can be expressed as
\begin{subequations}\label{eq:gen_dynamic_systems_sen_sel}
	\begin{align}
	\dot{\m x}(t) &= \mA \m x (t) + \m f(\m x) + \mB\m u(t),\\
	\m y(t) &= \m \Gamma(\m \gamma)\mC \m x (t).
	\end{align}
\end{subequations}
On top of that, it is also useful to define a set $\mathcal{G}\subseteq \{0,1\}^N$ to describe logistic constraints and availability of sensors so that we may impose $\m \gamma \in \mathcal{G}$.
With the addition of this constraint, a simplified high level formulation for SPP can be posed as 
\begin{align*}
\textbf{(P1)}\;\;\;	\minimize \;\;\;&  \m c^{\top} \boldsymbol \gamma + \mathrm{EstObjective} \\
\subjectto  \;\;\;& \eqref{eq:gen_dynamic_systems_sen_sel},\; \m \gamma\in \mathcal{G},\; \mathrm{EstConstraints}. 
\end{align*}
In~\textbf{P1} above, our objectives are threefold: \textit{(1)} performing state estimation for NDS \eqref{eq:gen_dynamic_systems_sen_sel} while \textit{(2)} utilizing smallest number of sensors as possible (or satisfying a given constraint over the collections of library of sensors) and \textit{(3)} optimizing a specific estimation metric (e.g., robust performance). The constant vector $\m c\in\mbb{R}^N_{+}$ in the objective function of \textbf{P1} assigns weights for each sensor $\gamma_i$. To achieve the above objectives, it is important to acquire some knowledge on the characteristics of nonlinearity of NDS. The ensuing section discusses approaches to parameterize $\m f(\cdot)$, which paves the way for a detailed SPP formulation. 

\setlength{\textfloatsep}{3pt}
\begin{table}[t]
		\vspace{-0.1cm}
	\centering
	\caption{Four classes of nonlinearity considered in this paper. 
	}\label{tab:nonlinear_class}
	\vspace{-0.2cm}
	\renewcommand{\arraystretch}{1.5}
	\begin{tabular}{|m{0.85in}|m{2.25in}|}
		\hline
		\textbf{Class} $\vphantom{\left(\frac{v}{l}\right)}$ & \textbf{Mathematical Property} $\vphantom{\left(\frac{v}{l}\right)}$ \\[0.5ex] \hline\hline
		\textit{Bounded Jacobian} &$\vphantom{\left(\frac{v_f^{G^0}}{l}\right)} -\infty < \barbelow{f}_{ij} \leq \dfrac{\partial f_i}{\partial x_j}(\m x)\leq \bar{f}_{ij} < + \infty $,
		\newline $\vphantom{\left(\frac{v}{l_j}\right)} \m x \in  \mathbfcal{X}$, $ \barbelow{f}_{ij}, \bar{f}_{ij} \in \mbb{R}$   \\ \hline
		\textit{Lipschitz} \newline \textit{Continuous} & $\vphantom{\left(\frac{v^i}{l}\right)}\norm{\m f(\m x)-\m f(\hat{\m x})}_2 \leq \beta \norm{\m x - \hat{\m x}}_2$, \newline  $\m x, \hat{\m x} \in \mathbfcal{X}$, $\beta \in \mbb{R}_{+}$ \\ \hline
		\textit{One-Sided Lipschitz} & $\vphantom{\left(\frac{v_f}{l}\right)} \left((\m f(\m x)-\m f(\hat{\m x}))\right)^{\top}\hspace{-0.1cm}(\m x-\hat{\m x})\leq \rho \norm{\m x - \hat{\m x}}_2^{2}$, \newline $\m x, \hat{\m x} \in  \mathbfcal{X}$, $\rho \in \mbb{R}$ \\ \hline
		\textit{Quadratically} \newline \textit{Inner-Bounded} & $\vphantom{\left(\frac{v_f}{l}\right)} \left((\m f(\m x)-\m f(\hat{\m x}))\right)^{\top}\hspace{-0.1cm}\left((\m f(\m x)-\m f(\hat{\m x}))\right) \leq $\newline $\delta_1 \norm{\m x - \hat{\m x}}_2^{2} +\delta_2\left((\m f(\m x)-\m f(\hat{\m x}))\right)^{\top}\hspace{-0.1cm}(\m x-\hat{\m x}) $,\newline $\m x, \hat{\m x} \in  \mathbfcal{X}$, $\delta_{1,2} \in \mbb{R}$ \\ \hline
	\end{tabular}
	\vspace{0.05cm}
\end{table}
\setlength{\floatsep}{3pt}

\section{Parameterizing Nonlinearity in NDS}\label{sec:NDS_classification}
As mentioned in the introduction, the nonlinearity of NDS \eqref{eq:gen_dynamic_systems} is assumed to belong to one of the following function sets: {bounded Jacobian}, {Lipschitz continuous}, {one-sided Lipschitz}, and {quadratically inner-bounded}---see Tab. \ref{tab:nonlinear_class} for the detailed mathematical definitions. If these conditions are satisfied inside the set $\mathbfcal{X}\subset \mbb{R}^{n_x}$, then we assert that the nonlinearity holds locally in $\mathbfcal{X}$. Otherwise the function sets are assumed to hold globally.
The premise that NDS \eqref{eq:gen_dynamic_systems} belongs to one of these function sets is actually not restrictive as it seems because of the following reasons. First, several major infrastructures, including power grids and traffic networks can actually be modeled as NDS which nonlinearities satisfy the conditions given in Tab. \ref{tab:nonlinear_class} \cite{Nugroho2018,Nugroho2019,Qi2018Access}. 
Second, advancements in observer design technique through linear matrix inequality (LMI) during past decades such as \cite{Phanomchoeng2010} for Lipschitz systems, \cite{zhang2012full} for one-sided Lipschitz and quadratically inner-bounded, and \cite{Jin2018} for bounded Jacobian systems, allow the state estimation problem for NDS with such nonlinearities to be solved efficiently via conventional solvers.

In this regard, we see that parameterizing nonlinear function $\m f(\cdot)$ is somehow crucial. 
The straightforward approach to parameterize/characterize $\m f(\cdot)$ is to just analytically compute the corresponding constants for each type of nonlinearity mentioned above. 
For instance, our previous works deal with analytical computation of Lipschitz constants for traffic networks \cite{Nugroho2018} and synchronous generator \cite{Nugroho2019}.
Nonetheless, for high-dimensional NDS having complex nonlinear dynamics, parameterizing the nonlinearities analytically can indeed be tedious and cumbersome. Even if the constant can successfully be computed, it could be too conservative.
As an alternative we may consider optimization-based methods to compute/approximate such constants given that the following assumption on the NDS holds.
\begin{asmp}\label{asmp:1}
	The NDS \eqref{eq:gen_dynamic_systems} satisfies following properties
	\vspace{-0.1cm}
	\begin{enumerate}
		\item $\mathbfcal{X} := \prod_{i=1}^{n_x} \mathbfcal{X}_i$ is nonempty where for each $i$, there exist $\barbelow{x}_i, \bar{x}_i\in\mathbb{R}$ such that $\mathbfcal{X}_i = \left[\barbelow{x}_i, \bar{x}_i\right]$. 
		\item  $\m f(\cdot)$ is differentiable and has continuous partial derivatives on $\mathbfcal{X}$.
	\end{enumerate}
\end{asmp}

The above assumptions allow us to parameterize $\m f(\cdot)$ in a much simpler way.
To show this, let $f_i:\mathbb{R}^{n_x}\rightarrow \mathbb{R}$ be a part of $\m f(\cdot)$ where its Lipschitz constant $\beta_i$ can be computed as
\begin{align}
\vspace*{-0.1cm}
\abs {f_i(\m x)-f_i(\hat{\m x})}_2 \leq \Big(\max_{\m z\in \mathbfcal{X}}\,\norm {\nabla f_i(\m z)}_2\Big)\norm {\m x-\hat{\m x}}_2.\label{eq:lip_grad}
\end{align}
Realize that this approach transforms the computation of $\beta_i$ to a global maximization problem as $\beta_i = \max_{\m z\in \mathbfcal{X}}\,\norm {\nabla f_i(\m z)}_2$. Once $\beta_i$ has been computed for all $i$, Lipschitz constant for $\m f(\cdot)$ can be determined by simply calculating $\beta = \sqrt{\sum^{n_x}_{i=1}\beta_i}$ \cite{Nugroho2019}.
In what follows we succinctly discuss two methods, referred to as \textit{stochastic point-based method} and \textit{interval-based method}, to compute $\beta_i$ as formulated in \eqref{eq:lip_grad}. Note that these methods are not limited for Lipschitz functions; they can also be utilized to parameterize other function sets listed in Tab. \ref{tab:nonlinear_class}.

\setlength{\textfloatsep}{3pt}
\begin{algorithm}[t]
	\caption{Numerical Approximation of $\beta_i$}\label{algo_LDS}
	\DontPrintSemicolon
	\textbf{input:} $\mathbfcal{X}$, $s$, $\norm {\nabla f_i(\cdot)}_2$\;
	\textbf{generate:} a sequence $\mathcal{S}(\m x, s)$ of LDS in $\mathbfcal{X}$ \;
	\textbf{initialize:}  $\beta_i \leftarrow  0 $\;
	\For{$j=1:s$}{
		$\beta_i \leftarrow \max(\beta_i, \norm {\nabla f_i(\m x_j)}_2)$,   
		$\m x_j \in \mathcal{S}(\m x, s)$ \;} 
	\textbf{output:} $\beta_i$
\end{algorithm}
\setlength{\floatsep}{3pt}

\vspace{-0.1cm}
\subsection{Stochastic Point-Based Optimization Approach}\label{ssec:sthochastic_opt}
In this section, we discuss a stochastic approach for solving global optimization problem posed in \eqref{eq:lip_grad}.
The basic idea here is to generate $s$ random points inside the set $\mathbfcal{X}$, which is then used to evaluate the corresponding objective function. 
This technique can be referred to as a \textit{Monte Carlo method}\cite{dalal2008}, which is in contrast with \textit{Quasi-Monte Carlo} methods that use a pseudo-random approach such as \textit{low-discrepancy sequences} (LDS). In essence, LDS are sequence of points that are distributed almost equally. As a consequence, any LDS has relatively small discrepancy, which in turn gives LDS a nice property over a pure random sampling method, explained as follows. First, define $f^*_i$ as the optimal value such that $f^*_i = \max_{\m z\in \mathbfcal{X}}\,\norm {\nabla f_i(\m z)}_2$. Then, suppose that there are $s$ number of points in $\mathbfcal{X}$ denoted by the sequence $ \mathcal{S}(\m x, s) := \left\{ \m x_j \right\}^s_{j = 1}$ where $\m x_j \in \mathbfcal{X}$. The approximation of $f^*_i$ is then given as $f^*_{i,s} := \max_{\m x\in \mathcal{S}(\m x, s)} \norm {\nabla f_i(\m x)}_2$. 
Since $\mathcal{S}(\m x, s)$ is a LDS, then  
$\lim_{s\to \infty} \mathit{D}^*(\mathcal{S}) = 0$,
where $\mathit{D}^*$ denotes the worst-case discrepancy \cite{dalal2008}. This yields $f^*_i = \lim_{s\to \infty} f^*_{i,s}$ \cite{kucherenko2005application}, suggesting that the value of $f^*_{i,s}$ is approaching the value of $f^*_i$ as more points of LDS are used. 
Some well known LDS are \textit{Halton}, \textit{Sobol}, and \textit{Niederreiter} sequences \cite{dalal2008}.
Algorithm~\ref{algo_LDS} describes a simple offline procedure to approximate $\beta_i$ via LDS. 

\setlength{\textfloatsep}{3pt}
\begin{algorithm}[t]
	\caption{Pseudocode of Interval-Based Algorithm}\label{pseudocode_IA}
	\DontPrintSemicolon
	\textbf{input:} $\mathbfcal{X}$, $\epsilon_h$, $ \epsilon_x$, $f^I_i(\cdot)$\;
	\textbf{initialize:} $\mathbfcal{S}_1 = \mathbfcal{X}$, $C =\{\mathbfcal{S}_1\}$, $l$, and $u$ based on $f^I_i(\cdot)$\;
	\textbf{set:} $\bar{\mathbfcal{S}}$ as a subset in $C$ with the greatest upper bound \;
	\textbf{\textit{Phase I}:} find the best upper bound \;
	\While{$u-l > \epsilon_h$ {\bf and}   $\abs{\bar{\mathbfcal{S}}}> \epsilon_x$}{
		Split $\bar{\mathbfcal{S}}$ via BnB routines\;
		Update $C$, $\bar{\mathbfcal{S}}$, $l$, and $u$ based on $f^I_i(\cdot)$\;
	}
	\textbf{\textit{Phase II}:} find the best lower bound\;
	\While{$u-l > \epsilon_h$ {\bf and}  $\exists \mathbfcal{S}_j\hspace{-0.00cm}\in C,\;\abs{\mathbfcal{S}_j}> \epsilon_x$}{
		Split ${\mathbfcal{S}}_j$ via BnB routines\;
		Update $C$, $\bar{\mathbfcal{S}}$, $l$, and $u$ based on $f^I_i(\cdot)$\;
	}
	\textbf{Output:} $\beta_i\leftarrow u$\;
\end{algorithm}
\setlength{\floatsep}{3pt}

\vspace{-0.1cm}
\subsection{Interval-Based Optimization Approach}\label{ssec:interval_opt}
Interval-based optimization utilizes the principle of \textit{interval arithmetic} to obtain the best (smallest) interval providing the upper ($u$) and lower ($l$) bounds for $f^*_i$. In this approach, global maximization problems are solved through \textit{branch and bound} (BnB) routines \cite{moa2007interval}, which are working based on the following principles.
In the branching step, the main problem is divided into smaller subproblems.
The corresponding upper and lower bounds of interval evaluation of each subsets are then computed accordingly. In the bounding step, several subsets that surely do not contain any maximizer is then discarded. These two routines are performed iteratively until the algorithm terminates.
In the context of Lipschitz constant, let $f^I_i(\cdot)$ be an interval extension of the objective function, i.e., $\norm {\nabla f_i(\cdot)}_2$---see \eqref{eq:lip_grad}, and $C$ be a \textit{cover} (a collection of subsets) in $\mathbfcal{X}$.
As we progress with BnB algorithm, we have to ensure that all maximizers are always contained in $C$ while $l$ and $u$ are updated based on the interval evaluation of $f^I_i(\cdot)$ for each subset of $C$. In this procedure, referred to as \textit{interval-based algorithm}, two constants $\epsilon_h$ and $ \epsilon_x$ are utilized: $\epsilon_h$ is useful to bound the interval containing the optimal value (if $u-l \leq \epsilon_h$, the algorithm terminates) while $\epsilon_x$ is useful to determine whether a certain subset $\mathbfcal{S}\in C$ can be split---if the maximum width of $\mathbfcal{S}$ is less than or equal to $\epsilon_x$, it will not be split further (as such, $\mathbfcal{S}$ is said to be \textit{atomic} \cite{moa2007interval}).
As the algorithm iterates, the distance between $l$ and $u$ is expected to be smaller and smaller. 
A high level pseudocode of the proposed interval-based algorithm is presented in Algorithm \ref{pseudocode_IA}. From this algorithm, Lipschitz constant $\beta_i$ is simply computed as $\beta_i \leftarrow u$. This approach, together with the stochastic point-based approach, can be used to parameterize the other function sets shown in Tab. \ref{tab:nonlinear_class} as well---the details of these approaches for NDS parameterization will be addressed in future work.
In the next section we discuss strategies on solving SPP given that the parameters of the nonlinearity have been determined.


\section{Convex MISDP Formulations of SPP for NDS}\label{sec:sensor_placement}
\vspace{-0.05cm}
\subsection{The Case for Lipschitz NDS}\label{ssec:sensor_placement_lipschitz}
After the NDS has been parameterized and a certain observer design has been chosen, we may proceed to address the sensor placement. In what follows, we show how to solve SPP given that the nonlinearity of the NDS is locally Lipschitz continuous in $\mathbfcal{X}$ with Lipschitz constant $\beta$ (note that the forthcoming strategies are \textit{not} limited for solving SPP for Lipschitz NDS---they can be extended to address SPPs for other types of NDS; see Section \ref{ssec:sensor_placement_general}). In particular, by considering an observer design approach introduced in \cite{Phanomchoeng2010}, the SPP can be posed as follows
\begin{subequations}\label{eq:sen_sel_obs}
			\vspace*{-0.05cm}
\begin{align}
&\hspace*{-0.2cm}\textbf{(P2)}\;	\minimize_{ \m P, \m Y, \m \gamma, \kappa}\;\;  \m c^{\top} \boldsymbol \gamma \\
&\hspace*{-0.2cm}\subjectto  \;\;	  \begin{bmatrix}
\m A ^{\top}\m P + \m P\m A   +\kappa\beta^2 \m I\\ - \m Y\m \Gamma (\m \gamma)\m C- \m C ^{\top}\m \Gamma (\m \gamma)\m Y ^{\top}& *\\
\m P & -\kappa \m I \end{bmatrix} \preceq 0 \label{eq:sen_sel_obs_1} \\ 
&\quad \quad \quad\quad\;\;\m P \succ 0, \;\kappa \geq 0,\; \m \gamma\in \mathcal{G},\;
 \m \gamma \in \{0,1\}^N.\label{eq:sen_sel_obs_2}\vspace*{-0.05cm}
\end{align}
\end{subequations}
In \textbf{P2} the objective is to minimize the number of the activated sensors while \textit{(a)} finding the observer gain matrix $\mL$ computed as $\mP^{-1}\mY$ and \textit{(b)} fulfilling the sensor logistic constraints. Since $\m \Gamma(\m \gamma)$ is a binary matrix variable, \textbf{P2} is a nonconvex optimization problem with \textit{mixed-integer bilinear matrix inequalities} due to the $\m Y \m \Gamma(\m \gamma)$ term. Solving \textbf{P2} returns the optimal combination of sensors $\m \gamma$ as well as the observer gain $\m L$ that guarantees the trajectory of estimation error dynamics 
\begin{align}
\dot{\m e}(t) = \left(\m A -\m L \m C\right)\m e(t)-\left({{\m f}}({\m x}(t))-{{\m f}}(\hat{\m x}(t))\right),\label{eq:est_err_dyn}
\end{align}
to converge towards zero. In \eqref{eq:est_err_dyn}, $\m e(t):= \m x (t)- \hat{\m x}(t)$ where $\hat{\m x}$ denotes the observer's ({estimated}) state. In what follows we will discuss our solution approach to solve \textbf{P2}.

In order to solve \textbf{P2} efficiently, it is crucial to identify the nonconvex terms that emerge in \eqref{eq:sen_sel_obs}. 
Notice that the nonconvexity of \textbf{P2} are twofold: first, the $0/1$ integer variables appearing in $\m\gamma$ (or $\m\Gamma(\m\gamma)$), and second, the multiplication of $\m\Gamma(\m\gamma)$ with $\m Y$. To solve \textbf{P2}, one reasonable approach is to transform \textbf{P2} into MISDP. The reformulation of \textbf{P2} from nonconvex MISDP to convex MISDP can be carried out by either using big-M method \cite{nugroho2018algorithms,Taha2018} or McCormick's relaxation \cite{chanekar2017optimal,mccormick1976computability}. 
With that in mind, here we present a way of reformulating \textbf{P2} into MISDP via McCormick's relaxation. 
This reformulation can be performed by defining a new matrix variable $\m Q := \m Y \m \Gamma(\m \gamma)$ where $\m Q\in\mbb{R}^{n_x\times n_y}$ and supposing that $\m Y$ (therefore, $\m Q$) are bounded such that $\barbelow{\m Y} \leq \m Y \leq \bar{\m Y}$. This allows \textbf{P2} to be posed as \textbf{P3}, described in the following problem
\begin{subequations}\label{eq:sen_sel_obs_new}
		\vspace*{-0.05cm}
	\begin{align}
	&\hspace*{-0.2cm}\textbf{(P3)}\;	\minimize_{ \m P, \m Y, \m Q, \m \gamma, \kappa}\;\;  \m c^{\top} \boldsymbol \gamma \\
	& \subjectto \;\; 	  \begin{bmatrix}
	\m A ^{\top}\m P + \m P\m A  \\ - \m Q\m C- \m C ^{\top}\m Q ^{\top}  +\kappa\beta^2 \m I& *\\
\m P & -\kappa \m I \end{bmatrix} \preceq 0 \label{eq:sen_sel_obs_new_1} \\ 
	 & \quad \quad \quad\quad\quad\;\m Q = \m Y \m \Gamma(\m \gamma), \;\barbelow{\m Y} \leq \m Y \leq \bar{\m Y},\; \eqref{eq:sen_sel_obs_2}.\label{eq:sen_sel_obs_new_2} 	\vspace*{-0.05cm}
	\end{align}
\end{subequations} 
Realize that as \textbf{P3} is more constrained than \textbf{P2}, then
 any solution of \textbf{P3} is always feasible for \textbf{P2}.
  Having formulated \textbf{P3}, in what follows, we present our main result which provides a convex MISDP reformulation towards \textbf{P3}.
  \vspace{-0.17cm}
\begin{theorem}\label{thm:misdp}
The nonconvex MISDP problem \textbf{P3} is equivalent to the following convex MISDP
\vspace{-0.05cm}
\begin{subequations}\label{eq:sen_sel_obs_misdp}
\begin{align}
	\textup{\textbf{(P4)}}\;	\minimize_{ \m P, \m Y, \m Q, \m \gamma, \kappa}\;\; &\m c^{\top} \boldsymbol \gamma \\[0.001\baselineskip]
	\subjectto \;\; &\eqref{eq:sen_sel_obs_new_1},\; \m P \succ 0,\;\m \gamma\in \mathcal{G},\;
	\m \gamma \in \{0,1\}^N \label{eq:sen_sel_obs_misdp_1}\\ 
&\m \Phi \,\m \xi(\m Q,\mY,\m \gamma, \kappa) \leq \m \nu, \label{eq:sen_sel_obs_misdp_2}
\end{align}
\end{subequations} 	
where $\m\Phi$, $\m\nu$, and $\m \xi(\m Q,\mY,\m \gamma, \kappa)$ are detailed in \eqref{eq:thm1_proof_final}. 
\end{theorem}
\vspace{-0.4cm}
\begin{proof}
Note that the constraint $\m Q = \m Y \m \Gamma(\m \gamma)$ in \textbf{P3} is equivalent to $Q_{(i,j)} = Y_{(i,j)}\gamma_j$ for all $i,j$ since $\m \Gamma(\m \gamma)$ is a block diagonal matrix that corresponds to $\m \gamma$. This yields the equivalency below
\begin{align}
Q_{(i,j)} = Y_{(i,j)}\gamma_j \Leftrightarrow 
Q_{(i,j)} = \begin{cases}
Y_{(i,j)}, &\text{if}\;\;\gamma_j = 1 \\
0, &\text{if}\;\;\gamma_j = 0, \label{eq:thm_1_proof_1}
\end{cases}\vspace{-0.0cm}
\end{align}
where $Y_{(i,j)}\in\left[\barbelow{Y}_{(i,j)},\bar{Y}_{(i,j)}\right]$ for all $i,j$. The transformation of the consequence in \eqref{eq:thm_1_proof_1} into convex MISDP is carried out through applying McCormick's relaxation, which can be explained as follows.
First, by realizing that the expressions $\bar{Y}_{(i,j)}-{Y}_{(i,j)}$, ${Y}_{(i,j)}-\barbelow{Y}_{(i,j)}$, $1-\gamma_j$, and $\gamma_j$ are all nonnegative and noting that $Q_{(i,j)} = Y_{(i,j)}\gamma_j$, we obtain the following series of inequalities
\begin{subequations}\label{eq:thm_1_proof_2} 
	\vspace{-0.4cm}
	\begin{align}
	\left(\bar{Y}_{(i,j)}-{Y}_{(i,j)}\right)\left(1-\gamma_j\right) &\geq 0  \nonumber \\
	\Leftrightarrow Q_{(i,j)} &\geq {Y}_{(i,j)}+\bar{Y}_{(i,j)}\left(\gamma_j-1\right) \label{eq:thm_1_proof_2a} \\
	\left(\bar{Y}_{(i,j)}-{Y}_{(i,j)}\right)\gamma_j &\geq 0  
\Leftrightarrow\bar{Y}_{(i,j)}\gamma_j \geq Q_{(i,j)}   \label{eq:thm_1_proof_2b} \\
\left({Y}_{(i,j)}-\barbelow{Y}_{(i,j)}\right)\left(1-\gamma_j\right) &\geq 0  \nonumber \\
\Leftrightarrow -Q_{(i,j)}   &\geq -{Y}_{(i,j)}+\barbelow{Y}_{(i,j)}\left(1-\gamma_j\right) \label{eq:thm_1_proof_2c} \\
	\left({Y}_{(i,j)}-\barbelow{Y}_{(i,j)}\right)\gamma_j &\geq 0  
\Leftrightarrow Q_{(i,j)} \geq \barbelow{Y}_{(i,j)}\gamma_j.   \label{eq:thm_1_proof_2d} 
	\end{align}
\end{subequations}
Second, by noticing that there are only two possible values of $\gamma_j$, substituting $\gamma_j = 1$ to \eqref{eq:thm_1_proof_2} yields
\begin{equation*}
\left.\begin{aligned}
{Y}_{(i,j)} \geq &\,Q_{(i,j)} \geq Y_{(i,j)}   \\
\bar{Y}_{(i,j)} \geq &\,Q_{(i,j)} \geq \barbelow{Y}_{(i,j)} 
\end{aligned}\right\} \Rightarrow Q_{(i,j)} = Y_{(i,j)}. \label{eq:eq:thm_1_proof_3}
\end{equation*}
On the other hand, substituting $\gamma_j = 0$ to \eqref{eq:thm_1_proof_2} yields
\begin{equation*}
\left.\begin{aligned}
0 \geq &\,Q_{(i,j)} \geq 0  \\
\bar{Y}_{(i,j)} \geq &\,Y_{(i,j)} \geq \barbelow{Y}_{(i,j)} 
\end{aligned}\right\} \Rightarrow Q_{(i,j)} = 0. \label{eq:eq:thm_1_proof_4}
\end{equation*}
Since this is true for all $i,j$, then as a result, \eqref{eq:thm_1_proof_1} and \eqref{eq:thm_1_proof_2} are equivalent.
Next define $\m\sigma_1$, $\m\sigma_2$, $\m \Psi$ and $\m \psi$ as follows
\begin{align*}
\m\sigma_1 &:= \bmat{1&-1&1&-1}^{\top},\;
\m\sigma_2 := \bmat{-1&1&0&0}^{\top} \\
\m \Psi &:= \m \Omega' \oslash \bmat{\m I_{n_y}\otimes\m 1_{4n_x}},\; \m\Omega':= \underbrace{\bmat{\mathrm{vec}(\m \Omega)&\hdots&\mathrm{vec}(\m \Omega)}}_{N\, \text{times}}   \\
\m \Omega &:= \hspace{-0.1cm}\bmat{\m \omega_{11}&\m \omega_{12}&\cdots & \m \omega_{1n_y}\\
	\m \omega_{21}&\m \omega_{22}&\cdots & \m \omega_{2n_y}\\
	\vdots&\vdots&\ddots&\vdots\\
	\m \omega_{n_x1}&\m \omega_{n_x2}&\cdots & \m \omega_{n_xn_y} }\hspace{-0.1cm},\,
\m \omega_{ij} := \hspace{-0.1cm}\bmat{-\barbelow{Y}_{(i,j)}\\\bar{Y}_{(i,j)}\\
	-\bar{Y}_{(i,j)}\\\barbelow{Y}_{(i,j)}}\hspace{-0.1cm},\,\forall i,j \\
\m \psi &:= \mathrm{vec}(\m \Omega)\oslash\left(\m 1_{n_xn_y}\otimes\bmat{1&1&0&0}^{\top}\right),
\end{align*}
where $\m I_{xy} := \m I_{n_xn_y}$. 
For all $i,j$, \eqref{eq:thm_1_proof_2} can be written as 
\begin{align*}
	&\left[
\begin{matrix}
\m I_{xy}\otimes \m\sigma_1
\end{matrix}  \,\vline  \,	
\begin{matrix}
\m I_{xy}\otimes \m\sigma_2 
\end{matrix} \,\vline\,  
\m \Psi
\right]\left[
\mathrm{vec}(\m Q)^{\top} \;\mathrm{vec}(\m Y)^{\top}\;\m \gamma^{\top} 
\right]^{\top}\leq \m \psi. 
\end{align*}
By combining the above with $\barbelow{\m Y} \leq \m Y \leq \bar{\m Y}$ and $\kappa \geq 0$, one can obtain \eqref{eq:sen_sel_obs_misdp_2} where $\m\Phi$, $\m\nu$, and $\m \xi(\m Q,\mY,\m \gamma, \kappa)$ are given as
\begin{subequations}\label{eq:thm1_proof_final}
\begin{align}
&\hspace{-0.3cm}\m\Phi := \mathrm{diag}\left(\left[
\begin{matrix}
\m I_{xy}\otimes \m\sigma_1
\end{matrix}  \,\vline  \,	
\begin{matrix}
\m I_{xy}\otimes \m\sigma_2 
\end{matrix} \,\vline\,  
\m \Psi
\right],\,\m I_{xy},\,-\m I_{xy},\,-1\right)\\
&\m\nu := \bmat{\m \psi^{\top} &\mathrm{vec}(\bar{\m Y})^{\top}&-\mathrm{vec}(\barbelow{\m Y})^{\top}& 0}^{\top}\\
&\m \xi(\m Q,\mY,\m \gamma, \kappa) := \bmat{\left[
	\mathrm{vec}(\m Q)^{\top} \;\mathrm{vec}(\m Y)^{\top}\;\m \gamma^{\top} 
	\right]^{\top} \\ \mathrm{vec}(\m Y)\\ \mathrm{vec}(\m Y)\\ \kappa}.
\end{align} 
\end{subequations}
This concludes the proof. 
\end{proof}
\vspace{-0.1cm}
Being a convex MISDP, \textbf{P4} can be conveniently solved using any general MISDP solver, such as BnB algorithm. Albeit it is known for its ability to return optimal solutions, unfortunately, its computational time, in general, does not scale well with the number of variables. This points out one major disadvantage of BnB algorithm. In this regard, the \textit{branch-and-cut} (BnC) algorithm can be explored to solve \textbf{P4}---this will be the focus of our future work. Essentially, BnC algorithm incorporates the advantages of BnB algorithm along with cutting plane method, thus providing a potentially faster alternative to solve convex MISDP problems.  

Other than the proposed approach, some other methods from the literature may also be considered. 
For instance, if the $0/1$ integer constraint is relaxed such that $\gamma_j \in  \left[0,1\right]$, \textbf{P2} can be solved through SDP routines---our earlier work \cite{Taha2018} presents a thorough study on addressing actuator selection for \textit{linear} dynamic systems through a robust control framework via SDP relaxations and approximations, as well as big-M method and heuristics. Recently, \cite{chang2018co} revisits the problem described in \cite{Taha2018} by proposing an approach that combines a new reformulation method to tackle the nonconvexity due to the multiplication between integer variables and matrix variable with BnB algorithm. By limiting the number of iterations on the BnB algorithm, the computational time can be shortened. Finally, an in-depth study that evaluates the performance of all these methods to solve the SPP for NDS can also be considered for future work.

\subsection{Extending The MISDP Formulation to Other Function Sets}\label{ssec:sensor_placement_general}
Having developed a convex MISDP to address the SPP for Lipschitz NDS, here we succinctly demonstrate a way to extend the proposed approach to solve similar problem for other types of nonlinearity as well as different observer structures. For the sake of illustration, consider the following observer design for bounded Jacobian systems developed using methods from \cite{Jin2018}
\begin{subequations}\label{eq:bj_obs}
	\vspace{-0.47cm}
	\begin{align}
	\mathrm{find}\;\;\; &\m P\succ 0, \;\m Y, \; \m \Lambda\geq 0\label{eq:bj_obs_1}\\
	\subjectto \;\;\; &\begin{bmatrix}
	\m A ^{\top}\m P + \m P\m A -\m Y\m C \\ - \m C ^{\top}\m Y ^{\top}  +\m \Theta_1\left(\m{\Lambda}\right)& * \\
	\m W^{\top}	\m P + \m\Theta_2\left(\m{\Lambda}\right)& \m\Theta_3\left(\m{\Lambda}\right)\end{bmatrix} \preceq 0,\label{eq:bj_obs_2}
	\end{align}
		where $\m\Theta_1\left(\m{\Lambda}\right)$, $\m\Theta_2\left(\m{\Lambda}\right)$, and $\m\Theta_3\left(\m{\Lambda}\right)$ are detailed as
	\begin{align*}
	\m\Theta_1\left(\m{\Lambda}\right) &=\hspace{-0.0cm} \mathrm{diag}\left(\left\{\sum_{i=1}^{n_x}\Lambda_{(i,j)}\left(\bar{c}_{ij}^2-\barbelow{c}_{ij}^2\right) \right\}^{n_x}_{j = 1} \right) \\
	\m\Theta_2\left(\m{\Lambda}\right)^{\top} \hspace{-0.1cm}&=\hspace{-0.0cm} \bmat{\left\{\mathrm{diag}\left(\left[\Lambda_{(i,1)}\barbelow{c}_{i1},\hdots,\Lambda_{(i,n_x)}\barbelow{c}_{in_x}\right]\right)\right\}^{n_x}_{i = 1}}\quad\quad\\
	\m\Theta_3\left(\m{\Lambda}\right) &=\hspace{-0.0cm} \mathrm{diag}\left(\mathrm{vec}\left(\m{\Lambda}\right)\right),
	\end{align*}
	where $\barbelow{c}_{ij}$ and $\bar{c}_{ij}$ are functions of the Jacobian bounds given as
	\begin{align*}
\barbelow{c}_{ij} = \frac{1}{2}\left(\barbelow{f}_{ij}+ \bar{f}_{ij}\right), \;\;
\bar{c}_{ij} = \frac{1}{2}\left(\barbelow{f}_{ij}- \bar{f}_{ij}\right),\;\; \forall i,j.
	\end{align*}
	\vspace{-0.0cm}
\end{subequations} 
In \eqref{eq:bj_obs_1}, the constraint $\m \Lambda \geq 0$ indicates that each element of $\m \Lambda$ has to be nonnnegative, i.e., $\Lambda_{(i,j)} \geq 0$.
In the context of SPP with integer variable $\m \gamma$, then due to the sensor selection, the term $\m Y \m C$ in \eqref{eq:bj_obs_2} becomes $\m Y\m \Gamma (\m \gamma)\m C$. By again supposing that the variable $\m Y$ is bounded such that $\barbelow{\m Y} \leq \m Y \leq \bar{\m Y}$ and defining a new matrix variable $\m Q := \m Y \m \Gamma(\m \gamma)$, we can utilize the same trick presented in Theorem \ref{thm:misdp} to transform this problem into convex MISDP.  
The next section showcases the effectiveness of convex MISDP posed in \textbf{P4} to solve SPP via BnB routines. 

\begin{figure}[t]
	\centering
	\vspace{-0.3cm}
{\includegraphics[scale=0.53]{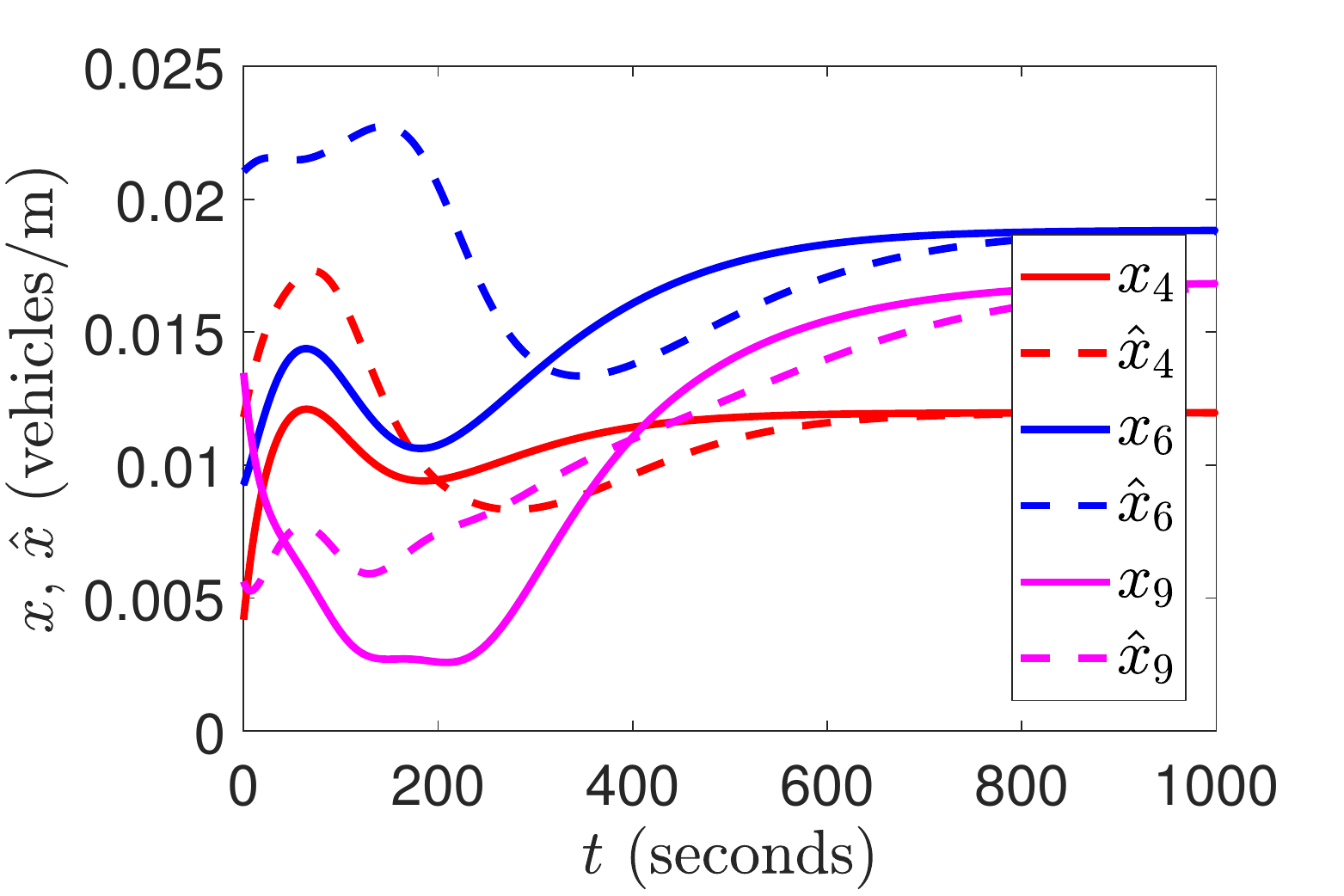}\label{Fig:state_1}}{\vspace{-0.05cm}}\hspace{-0.9em}
	\vspace{-0.15cm}
	\caption{Dynamic state estimation results for the traffic network assuming free-flow condition: a comparison between system's (actual) states $\m x$ and observer's (estimated) states $\hat{\m x}$ for the $\nth{4}$, $\nth{6}$, and $\nth{9}$ highway segments.}
	\label{Fig:result}	
	\vspace{-0.0cm}
\end{figure}

\section{Numerical Example: SPP on Traffic Networks}\label{sec:numerical}
In this section, we illustrate the proposed method for solving SPP for traffic networks on stretched highways. The model of traffic networks on a stretched highway consisting inflow and outflow ramps can be constructed by creating some partitions called \textit{segments} on the highway of equal length $l$ such that the dynamics of each segment can be expressed as \cite{Contreras2016}
\begin{align*}
\dot{\rho}_i (t)  
&=  \frac{1}{l}\left(\sum q_{i-1}(t)-\sum q_i (t)\right),
\end{align*}
where $q_{i-1}$, $q_{i}$ represent any inflow and outflow and ${\rho}_i$ represents traffic density on that segment. This model assumes that the highway is in free-flow condition. As such the traffic density satisfies $\rho_i\in[0,\rho_c]$, where $\rho_c$ denotes the critical density. 
A complete NDS model of this traffic network, which is proved to be locally Lipschitz continuous, is available in \cite{Nugroho2018} and can be expressed in the form of \eqref{eq:gen_dynamic_systems}.  
In this model, the state $\m x$ represents the traffic density on all segments, where the inflows and outflows are perceived as input $\m u$. The output $\m y$ is chosen to be equal to state $\m x$ such that $\m C =\m I$ since traffic sensors are considered measuring traffic density on each highway segment. The nonlinearities in $\m f(\cdot)$ takes a quadratic form---readers are referred to \cite{Nugroho2018} for the detailed parameters and expressions. 

The aim of this numerical test is to seek the minimum configuration of traffic sensor such that a dynamic state estimation on all highway segments can be performed.
This particular highway is assumed to consist $N_M=10$ highway segments, $N_I=2$ on-ramps, and $N_O=4$ off-ramps. As such the total number of states are equal to $n_x = N = N_M + N_I+N_O$.
The following parameters, adapted from \cite{Contreras2016}, are considered: $v_f = 31.3$ m/s, $\rho_m = 0.053$ vehicles/m, and $l = 500$ m. The exit ratio is chosen to be equal to $\alpha = (0.2,0.3,0.4,0.5)$ while the steady-state flow vector is set to be $u = [0.2\;0.1\hspace{-0.05cm}\times\hspace{-0.05cm}\m 1_{1\times 2}\;0.01\hspace{-0.05cm}\times\hspace{-0.05cm}\m 1_{1\times 4}]^{\top}$.
The values of $\barbelow{\m Y} $ and $\bar{\m Y} $ in \textbf{P4} are respectively chosen to be $-10^2\times\m 1_{n_x\times n_y}$ and $10^2\times\m 1_{n_x\times n_y}$, whereas the weighting vector $\m c$ is set to be $\m 1$. To simulate a more realistic condition, we impose a logistic constraint such that $1 \leq \sum_i \gamma_i \leq 8$.
The Lipschitz constant for this case is analytically computed using a formula given in \cite{Nugroho2018}, which value is equal to $0.34510$. The computed Lipschitz constant using interval-based optimization is equal to $0.29362$, which is less conservative than the analytical Lipschitz constant. In this numerical example, we opt to implement the analytical Lipschitz constant.
The initial conditions for the dynamic simulation are randomly generated such that $0 \leq x_i(0) \leq \rho_c$. The simulations are performed using MATLAB R2017b running on a 64-bit Windows 10 with 2.5GHz Intel\textsuperscript{R} Core\textsuperscript{TM} i7-6500U CPU and 16 GB of RAM.
YALMIP's \cite{Lofberg2004} BnB algorithm along with MOSEK's \cite{andersen2000mosek} SDP solver
are utilized to solve \textbf{P4}. 

After \textbf{P4} is successfully solved, the obtained optimal sensor placement is $\m\gamma^* = [\m 0_{1\times 15} \;1]^{\top}$.
The BnB algorithm terminates after 2 iterations with total computation time of $1.58$ seconds.
The corresponding observer gain $\m L$ together with $\m\gamma^*$ produce an asymptotically stable estimation error dynamics as depicted in Fig.~\ref{Fig:result}. This means that only 1 out of 8 available sensors is needed by the observer to perform traffic density estimation.   
It is important to note that, according to the theorem presented in \cite{Contreras2016}, in a free-flow condition, there should be one sensor measuring the last (downstream) segment to achieve an observable system.
Therefore our result corroborates this theorem, as we only use one sensor to perform dynamic estimation. 
However, notice that the result in \cite{Contreras2016} is intended for a \textit{linearized} traffic network model, in contrast to nonlinear model considered here.
With that in mind, our work provides solutions for some issues in this particular research direction that were not addressed in previous literature, thus showcasing an additional contribution made in this work.

%
\section{Summary and Future Work}\label{sec:conclusion_future_work}
\vspace{-0.00cm}
This paper discusses some strategies to address the SPP for NDS which nonlinearity belongs to the following function sets: bounded Jacobian, Lipschitz continuous, one-sided Lipschitz, and quadratically inner-bounded.
We demonstrate that the SPP can be solved through a hybrid of Lyapunov theory and mixed integer convex programming. Specifically, we aim to \textit{(a)} find the least number of sensors (along with the location of sensors) while \textit{(b)} finding the corresponding observer gain matrix. A numerical example on a simple traffic network showcases the effectiveness of the proposed approach. Future work will include developing and implementing a special type of BnC algorithm to provide a potentially more scalable alternative, compared to BnB algorithm, to solve SPP especially for large-scale NDS, in addition to incorporating uncertainty and robustness metrics.

\vspace{-0.0cm}

\bibliographystyle{IEEEtran}	\bibliography{bib_file}

\end{document}